\documentclass[10pt,twoside]{article}
\usepackage{graphicx}
\usepackage{amsmath}
\usepackage{Latex-document}

\title{\bf Bubbling and Regularity Issues  \vskip -2mm
in Geometric Non-linear Analysis \vskip 6mm}

\author{T. Rivi\`ere\vspace*{-0.5cm}\thanks{D-Math., ETH-Zentrum,
CH-8092 Z\"urich, Switzerland. E-mail: riviere@math.ethz.ch}}

\date{\vspace{-8mm}}

\newcommand{\be}{\begin{equation}}
\newcommand{\ee}{\end{equation}}
\def\R{{{{\rm l} \kern -.15em {\rm R}}}}
\def\Z{{{{\rm Z} \kern -.35em {\rm Z}}}}

\def\lf{\left}
\def\rg{\right}

\def\ep{\varepsilon}

\def\ov{\overline}
\def\Om{\Omega}
\def\om{\omega}
\def\p{\partial}

\begin{document}

\maketitle

\thispagestyle{first} \setcounter{page}{197}

\begin{abstract}\vskip 3mm

 Numerous elliptic and parabolic variational problems arising
 in physics and geometry (Ginzburg-Landau equations, harmonic
  maps, Yang-Mills fields, Omega-instantons, Yamabe equations,
 geometric flows in general...) possess a critical dimension in
  which an invariance group (similitudes, conformal groups) acts.
 This common feature generates, in all these different situations,
 the same non-linear effect. One observes a strict splitting in
 space between an almost linear regime and a dominantly non-linear regime
 which has two major characteristics : it requires a quantized amount of energy
and
 arises along rectifiable objects of special geometric
 interest (geodesics, minimal surfaces, J-holomorphic curves,
 special lagrangian manifolds, mean-curvature flows...).

\vskip 4.5mm

\noindent {\bf 2000 Mathematics Subject Classification:} 35D10,  35J20, 35J60, 49Q20, 58E15, 58E20.

\noindent {\bf Keywords and Phrases:} Harmonic maps, Yang-Mills fields, Ginzburg-Landau vortices.
\end{abstract}

\vskip 12mm

\section{Energy quantization phenomena for harmonic maps and Yang-Mills Fields}
 \label{section 1}
\setzero\vskip-5mm \hspace{5mm}

\subsection{The archetype of energy quantization: the {\boldmath $\ep$}-regularity}

\vskip-5mm \hspace{5mm}

Let $B^m$ be the flat $m$-dimensional ball and $N^n$ be a compact, without boundary Riemannian manifold.
By the Nash embedding Theorem, we may assume  $N\subset{\R}^k$. Let $W^{1,2}(B^m,N^n)$ be the maps in $W^{1,2}(B^m,{\R}^k$ that take values
 almost everywhere in $N^n$. A map $u\in W^{1,2}(B^m,N^n)$ is called stationary harmonic if it is a critical point
 for the Dirichlet energy
\[
E(u)=\int_{B^m}|\nabla u|^2\ dx
\]
for both perturbations in the target (of the form $\pi_N(u+t\phi)$ for any $\phi$ in $C^\infty_0(B^m,{\R}^k)$
where $\pi_N(y)$ is the nearest neighbor of $y$ in $N$)  and in the domain (of the form $u\circ(id+t X)$ for any
vector field $X$ in $C^\infty_0(B^m,{\R}^m)$).
There are relations between these two conditions, in particular a smooth weakly harmonic map $u$
(i.e. critical for perturbations in the target) is automatically stationary (i.e. critical for perturbations in the domain).
 This is not true in the general
case : there exists weakly harmonic maps which are not stationary (see \cite{HLP}) and that can even be nowhere
continuous
 (see \cite{Ri1}).
 On the one hand, as a consequence of
being weakly harmonic, one has the Euler Lagrange equation ({\it harmonic map equation})
\be
\label{3}
\Delta u+A(u)(\nabla u, \nabla u)=0\quad ,
\ee
where $A(u)$ is the second fundamental form of $N$ embedded in ${\R}^k$.
Stationarity, on the other hand, implies the {\it monotonicity formula} saying that for any
point $x_0$ in $B^m$ the {\it density of energy} ${r^{2-m}}\int_{B_r(x_0)}|\nabla u|^2$ is an increasing function
of $r$. Among stationary harmonic maps are the {\it minimizing harmonic maps} (minimizing $E$ for their boundary datas).
In \cite{ScU} R. Schoen and K. Uhlenbeck proved  that
there exists $\ep(m,N)>0$ such that, for any minimizing harmonic map, $u$ in $W^{1,2}(B^m,N^n)$ and any ball
$B_r(x_0)\subset B^m$ the following holds
\[
\frac{1}{r^{n-2}}\int_{B_r(x_0)}|\nabla u|^2\ dx\le \ep(m,N)
\quad\quad\Longrightarrow \quad\quad \|\nabla
u\|_{L^\infty(B_{r/2}(x_0))}\le \frac{C}{r}
\]
where $C$ only depends on $m$ and $N^n$. In other words, there exists a
number $\ep>0,$ depending only $m$ and $N$, making a strict splitting
between the {\it almost linear regime} for minimizing harmonic maps in
which derivatives of the maps are under control and the {\it totally
non-linear regime} where the map may be singular. Since this theorem in
the early eighties, and previous partial regularity result of that type
for minimal surfaces and other elliptic systems by Almgreen, DeGiorgi,
Giusti, Miranda, Morrey, $\cdots$, $\epsilon-$regularity theorems has
been found in various problems of geometric analysis ( geometric flows,
Yang-Mills Fields, $\cdots$, etc). In particular it was a natural
question whether the above result was true for arbitrary stationary
harmonic maps. This was done 10 years later by L.C. Evans (\cite{Ev}) in
the case where $N^n$ is a standard sphere. Evans' result was extended to
general target by F. Bethuel in \cite{Be}.
 Evans benefitted from a deeper understanding of the non-linearity in (\ref{3})
 developed in H\'elein's proof that any weakly harmonic map
from a two dimensional domain is $C^\infty$ ( see \cite{He} and \cite{CLMS}). This fact will be discussed in the next subsection.

Combining the $\ep-$regularity result with a classical Federer-Ziemer covering argument yields the following
upper-bound on the size of the singular set of stationary harmonic maps.
For a stationary harmonic map $u$ let Sing$\,u$  be the complement of the largest
 open subset where $u$ is $C^\infty$. Then
\be
\label{5}
{\mathcal H}^{m-2}(\mbox{Sing}\,u)=0\quad ,
\ee
where ${\mathcal H}^{m-2}$ denotes the $m-2-$dimensional Hausdorff measure. For stationary maps and general targets,
this is the best estimate available.
For minimizing harmonic maps Schoen and Uhlenbeck proved the optimal result :
dim(Sing$\, u)\le m-3$.
The reason for this improvement  is the striking fact that weakly converging
{\it minimizing harmonic maps} in $W^{1,2}(B^m,N)$ are in fact strongly converging in this space. This feature is
very specific to the minimizing map case. When we blow up the map at a singular point strong convergence then implies roughly
that the ``singular set of the limit is not smaller than the limit of the singular set''.{\it The Federer dimension reduction argument}
from the theory of minimal surfaces then gives by induction the result.
Here one sees the strong connection between understanding the singular set  and compactness properties
of sequences of solutions and  bubblings as presented in the next subsection.

\subsection{Bubblings of harmonic maps}
\vskip -5mm \hspace{5mm}

The following result of J.Sacks and K.Uhlenbeck \cite{SaU} is perhaps the earliest example of energy quantization
 in non-linear analysis. Their aim was to extend the work of Eells and Sampson \cite{ES} to
targets of not necessarily negative sectional curvatures, that is : find in a given 2-homotopy class
of an arbitrary riemannian manifold $N^n$ a more ``natural'' representant that would minimize
 the area and the Dirichlet energy $E$ in the given homotopy class . This raises the question of compactness
and the possible reasons for the lack of compactness for  harmonic maps from a 2-sphere (or
even more generally a Riemann surface) into $N^n$. We should mention not only the
 original contribution of J.Sacks and K.Uhlenbeck that focused on minimizing sequences but also
the later works by J.Jost \cite{Jo} for critical points in general, by M. Struwe \cite{St} for it's heat flow version
and the more recent contributions by T.Parker \cite{Pa}, W. Ding and G. Tian \cite{DiT} and F.H.Lin and
C. Wang \cite{LW}. The following result has influenced deeply the
non-linear analysis of the eighties, from the concentration-compactness of P.L.Lions to the
compactness of J-holomorphic curves by M.Gromov and the analysis of self-dual instantons on 4-manifolds
by S.K. Donaldson and K.Uhlenbeck.

{\bf Theorem 1.} \cite{SaU}, \cite{Jo}
\it Let $u_n$ be a sequence of weakly harmonic maps from a surface $\Sigma$ into a closed manifold
$N^n$ having a uniformly bounded energy. Then  a subsequence $u_{n'}$ weakly
converges in $W^{1,2}(\Sigma,N^n)$ to a harmonic map $u$ into $N^n$. Moreover, there exist
finitely many points $\{a_1\cdots a_k\}$ in $\Sigma$ such that the convergence
is strong in $\Sigma\setminus\{a_1\cdots a_k\}$ and the following holds
\[
|\nabla u_{n'}|^2\ dvol_{\Sigma}\rightharpoonup |\nabla u|^2\ dvol_{\Sigma}+\sum_{j=1}^k m_j\delta_{a_j}\quad
\quad\mbox{in Radon measure}
\]
where $m_j=\sum_{i=1}^{P_j}E(\phi^j_i)$ and $\phi^j_i$ are nonconstant harmonic 2-spheres of $N^n$
(harmonic maps from $S^2$ into $N^n$).\rm

The loss of energy during the weak convergence is not only concentrated at points but is also quantized :
the amount is given by a sum of energies of harmonic 2-spheres of $N^n$, the {\it bubbles},
 that might sometimes be even explicitly known
 (for instance if $N=S^2$, $E(\phi^j_i)\in 8\pi{\Z}$). The striking fact in this result is that it excludes
the possibility of losing energy in the neck between $u$ and the bubbles or between the bubbles themselves.
 This {\it energy identity} is quite surprising. It is a-priori conceivable, for instance, that, in a tiny annulus
 surrounding a blow-up point, an axially symmetric harmonic map
into $N$, that is a portion of geodesic of arbitrarily small length in
$N$, breaking the quantization of the energy, would appear. The {\it
energy identity} disappears if instead of exact solutions we consider in
general  Palais-Smale sequences for $E$ in general (see the work of
Parker [Pa] that was following his bubbling picture he established in
collaboration with J.Wolfson for pseudo-holomorphic curves).

Only relatively recently a first breakthrough was made by F.H. Lin in the attempt of extending Sacks Uhlenbeck result
 beyond the conformal dimension.

{\bf Theorem 2.} \cite{Li}
\it Let $u_n$ be a sequence of stationary harmonic maps from $B^m$ into a closed manifold
$N^n$ having a uniformly bounded energy, then there exists a subsequence $u_{n'}$ weakly
converging in $W^{1,2}(B^m,N^n)$ to a map $u$ and there exists a $m-2$ rectifiable subset
$K$ of $B^m$ such that the convergence is strong in $B^m\setminus K$ and moreover
\[
|\nabla u_{n'}|^2\ dvol_{\Sigma}\rightharpoonup |\nabla u|^2\ dvol_{\Sigma}+ f(x)\ {\mathcal H}^{m-2}\lfloor K\quad
\quad\mbox{in Radon measure}
\]
where $f$ is a measurable positive function of $K$.\rm

This result, establishing the regularity of the blow-up set of weakly converging stationary harmonic map,
is related to the resolution by D.Preiss \cite{Pr} of the Besicovitch conjecture on measures admitting densities.
Given a positive Borel regular measure $\mu$ such that there exists an integer $k$ for which $\mu$-a.e. the density
$\lim_{r\rightarrow 0}\mu(B_r)/r^k$ exists and is positive, it is proved in \cite{Pr} that there is a rectifiable
$k-$dimensional rectifiable set $K$ such that $\mu(B^m\setminus K)=0$. In the present situation calling
$\mu_n=|\nabla u_{n}|^2\ dvol_{\Sigma}$ converging to the Radon measure $\mu=|\nabla u|^2\ dvol_{\Sigma}+\nu$,
from the monotonicity formula one deduces easily that the defect measure $\nu$ fulfills the assumptions
of Preiss theorem for $k=m-2$ and the rectifiability follows at once. It has to be noted that the original proof
of Theorem 2 in \cite{Li} is self-contained and avoids D. Preiss above result.

Beyond the regularity of the defect measure $\nu$, the question  remained whether the whole picture established
 in the conformal dimension could be extended to higher dimensions and if  the {\it energy identity} still
holds. F.H. Lin and the author brought the following answer to that question.

{\bf Theorem 3.} \cite{LR2}
\it
Let $f$ be the function in the previous theorem, if $N^n=(S^n,g_{stand})$ or
if there is a uniform bound on $\|u_{n'}\|_{W^{2,1}(B^m)}$, then, for ${\mathcal H}^{m-2}$
a.e. $x$ of $K$, there exists a finite family of harmonic 2-spheres in $N^n$
$(\phi^j_x)_{j=1\cdots P_x}$ such that
\[
f(x)=\sum_{j=1}^{P_x} E(\phi^j_x).
\]
\rm

The proof of this {\it energy identity} in higher dimension is of
different nature from the one provided previously in dimension 2 where
the use of objects relevant to the conformal dimension only, such as the
Hopf differential, was essential. The idea in higher dimension was to
develop a technique of slicing, averaging method combined with estimates
in Lorentz spaces $L^{2,\infty}-L^{2,1}$. This technique seems to be
quite general as it have some genericity and permitted to solve problems
of apparently different nature as we shall expose in the next section.
The requirement of the $W^{2,1}$ bound
 for the case of a general target seems to be technical and should be removed.
 In the particular case when the target is the round sphere it was
proved in \cite{He} and \cite{CLMS} that the non-linearity $A(u)(\nabla u, \nabla u)$
in equation (\ref{3}) is in the Hardy space ${\mathcal H}^1_{loc}$ which immediately implies
the desired $W^{2,1}$ bound for the maps $u$. Whether this fact and this
$W^{2,1}$ bound can be extended to general targets is still unknown.

 The 2 previous results  suggest to view the loss of
compactness in higher dimension as being exactly the one happening in the conformal $2$-dimensional
case in the plane normal to $K$, locally invariant in the remaining $m-2$ dimensions tangent to $K$.
This understanding of the loss of compactness through creation of 2-bubbles for stationary
harmonic maps has consequences in regularity theory. Indeed, in the case where $N^n$
admits no harmonic 2-spheres (take for instance $N^n$ a surface of positive genus),
 no blow-up can arise ; weakly converging stationary harmonic maps are
strongly converging and the {\it Federer dimension reduction argument} discussed
at the end of the previous subsection combined with the analysis in \cite{Si} may be applied to improve the bound of
the size of the singular set  to dim$\,$Sing$\, u\le m-4$ (see \cite{Li}).

\subsection{High dimensional gauge theory}
\label{section 3}
\vskip -5mm \hspace{5mm}

The work of S.K. Donaldson and R. Thomas \cite{DoT} has given a new boost in the motivation for developing
the non-linear analysis of high dimensional gauge theory.

\noindent{\bf Bubbling of Yang-Mills Fields}

 Theorems 1,2,3 and their proofs are transposable to many others
 geometric variational problems
(such has Yang-Mills fields, Yamabe metrics...) having a given conformal invariant dimension $p$
(p=2 for harmonic maps, 4 for YM...etc). For instance in \cite{Ti}, G.Tian established the result corresponding
to theorem 2 for Yang-Mills fields.

 Consider a vector bundle $E$ over
a Riemannian manifold $(M^m,g)$ and assume $E$ is issued from a principal bundle
whose structure group is a compact Lie group $G$. Yang-Mills fields are connections $A$
on $E$ whose curvature $F_A$ solves
\be
\label{I.1}
d^\ast_A F_A=0
\ee
where $d^\ast_A$ is the adjoint to the operator $d_A$, acting
on $\wedge^\ast M\otimes End\, E$, with respect to the metric $g$ on $M$ and the Killing form of $G$ on the
fibers (recall that the Bianchi identity reads instead $d_A F_A=0$). Yang-Mills fields are
the critical points  of the Yang-Mills functional $\int |F_A|^2\ dvol_g$ for perturbations of the form $A+ta$ where
 $a\in \Gamma(\wedge^1M\otimes ad\,Q)$. It is proved in \cite{Ti} that, taking a sequence of smooth
Yang Mills connections $F_{A_n}$ having a uniformly bounded Yang-Mills energy, one may extract a subsequence
(still denoted $A_n$) such that for some $m-4$ rectifiable closed subset $K$ of $M$ the following holds :
 in a neighborhood of any point
of $M\setminus K$ there exist good choices of gauges such that $A_n$,
expressed in these gauges,  converges in $C^k$ topology (for any $k$) to
a limiting form $A$ that defines globally a smooth Yang-Mills connection
on $i^{\ast}E$, where $i$ is the canonical embedding of $M\setminus K$
in $M$. Moreover, one has the following convergence in Radon measure :
\be \label{I.2} |F_{A_n}|^2\ dvol_{M}\rightharpoonup |F_A|^2\ dvol_{M}+
f(x)\ {\mathcal H}^{m-4}\lfloor K\quad , \ee where $f$ is some
non-negative measurable function on $K$. For $m=4$, or under the
assumption that the measures $|\nabla_{A_n}\nabla_{A_n} F_{A_n}|\
dvol_{M}$ remain uniformly bounded it is proved by the author in
\cite{Ri4} that $f(x)$ is, ${\mathcal H}^{m-4}$ almost everywhere on
$K$, quantized and equal to a sum of energies of Yang-Mills connections
over $S^4$ : the {\it energy identity} holds. The proof use again the
Lorentz space duality $L^{2,\infty}-L^{2,1}$ applied this time to the
curvature. These estimates are consequences of $L^{4,2}$ estimates on
the connections $A$ for $m=4$. It is interresting to observe that such
an $L^{4,2}$ estimate of the connection plays a crucial role in the
resolution by J.Shatah and M. Struwe of the wave map Cauchy problem in 4
dimension
 for small initial datas in the ``natural'' space $H^2\times H^1$ (see \cite{ShS}).

To complete the description of the blow-up phenomena several open questions remain both for
harmonic maps and Yang-Mills fields. First, what are the exact nature of the limiting map $u$ and
connection $A$ ? Is $u$ still a stationary harmonic map on the whole ball $B^m$ ? Is  $A$ still a
smooth Yang-Mills connection of some new bundle $E_0$ over $(M,g)$ ?
Can one expect in both cases more regularity than the rectifiability for $K$ ?
What is the exact nature of the concentration set $K$ ?

It happens that these two questions are strongly related to eachother :
For instance, in the context of harmonic maps, if the weak limiting map
$u$ was also a stationary harmonic map, then one would deduce from the
monotonicity formula that $K$ is a stationary varifold and therefore
inherits nice regularity properties. There is a similar notion of
stationary Yang-Mills field (see \cite{Ti}). The limiting $A$ being
stationary, as expected, likewise would imply the stationarity of $K$.
These questions in both cases are still widely open and seem to be
difficult. Nevertheless in the case of Yang-Mills, in a work in
preparation, T. Tao and G. Tian are proving a {\it singularity
removability} result saying that the weak limit $A$ of \underbar{smooth}
Yang-Mills $A_n$ can be extended to a  smooth Yang-Mills connection
aside from a ${\mathcal H}^{m-4}$ measure zero set. Their proof is based
on an $\epsilon$-regularity result for weak Yang-Mills in high dimension
which was also independently obtained by Y. Meyer and the author in
[MR].

An important case where these questions have been solved is :

\noindent{\bf The case of {\boldmath $\Omega$}-anti-self-dual
Instantons}

This notion extends the 4-dimensional notion of instantons to
higher dimension. Assuming $m\ge 4$, and given a closed
$m-4$-form $\Omega$, we say that a connection $A$ is an
$\Omega-$anti-self-dual connection when the following equation
holds \be \label{I.3} -\ast(F_A\wedge\Omega)=F_A. \ee
In view of
the Bianchi identity $d_A F_A=0$ and the closedness of $\Omega$
it follows that $d^\ast_AF_A=0$ that is $A$ is a particular
solution to Yang-Mills equations. It is then shown in \cite{Ti}
that, if we further assume that the co-mass of $\Omega$ is less
than 1, then $K$ in (\ref{I.2}) is a minimizing current in
$(M,g)$ calibrated by $\Om$ (i.e. $\Om$ restricted to $K$
coincides with the volume form on $K$ induced by $g$) and
therefore minimizes the area in it's homology class.

A special case of interest arises when $(M,g,\om)$ is a Calabi-Yau 4-fold
with a global ``holomorphic volume form'' $\theta$ (i.e. $\theta$ is a holomorphic $({4}, 0)-$form satisfying
$\theta\wedge\ov{\theta}=dvol_{M}$). Let $\Om$ be the form generating the $SU(4)$ holonomy on $M$
(i.e. $\Om$ is a unit holomorphic section of the canonical bundle $\wedge^{(4,0)}T^\ast M$) and assume $M$ is the product
$T\times V^3$ of a 2-torus $T$ with a Calabi Yau 3-fold $V^3$, for a given $SU(2)-$bundle E, $T-$invariant solutions to (\ref{I.3}) are  pairs
 of connections $B$ on $E$ over $V^3$ and sections $\phi$ of $End\, E$  that solve a
vortex equations (see 6.2.2 in \cite{Ti}). The moduli space of these solutions should be related to the so-called
 holomorphic Casson invariants (see also R.Thomas work \cite{Th}). The loss of compactness of these solutions arises along $T$ times holomorphic curves
 in $V^3$.
To relate this moduli space with the holomorphic Casson invariants it is necessary to perturb, in a generic way,
 the complex structure into a not necessarily integrable one. This generates several non-linear analysis questions
 related to bubbling and regularity which have to be solved :
\begin{itemize}
\item[i)] Show that the weak limits of the vortex equations in the almost complex setting in $V^3$
have only isolated point singularities located on the pseudo-holomorphic blow-up set
(see the extended conjectures on the singular sets of solutions to (\ref{I.3}) in \cite{Ti}).
\item[ii)] Construct solutions to the vortex equations which concentrate on some
given choice of pseudo-holomorphic curves.
\item[iii)] Show that an arbitrary 1-1 rectifiable cycle (i.e. a 2 dimensional cycle whose tangent plane is invariant
 under the almost complex structure action) in an almost complex manifold (that may arise
as a blow-up set in the above case) is a smooth surface aside from eventually isolated branched points.
\end{itemize}

In collaboration with G.Tian, in the direction of iii), the following
result was established.

{\bf Theorem } \cite{RT} \it Let $(M^4,J)$ be a smooth almost complex
4-real manifold. Then any 1-1 rectifiable cycle is a smooth
pseudo-holomorphic curve aside from isolated branched points. \rm

\section{Ginzburg-Landau line vortices}
\label{section 2} \setzero\vskip-5mm \hspace{5mm}

\subsection{The strongly repulsive asymptotic}

\vskip-5mm \hspace{5mm}

The free Ginzburg-Landau energy on the 3-dimensional ball reads \be \label{6}
GL(u,A)=\int_{B^3}|d\,u-iAu|^2+\frac{\kappa^2}{2}(1-|u|^2)^2+|dA|^2 \ee where $u$ is a complex function satisfying
$|u|\le 1$ and $A$ is a 1-form on $B^3$. The study of the above variational problem was initiated by A. Jaffe and
C. Taubes in 2 dimensions and by J. Fr\"ohlich and M. Struwe in \cite{FS} in 3 dimensions. The connected
components of the zero set of the order parameter $u$ are the so called {\it vortices} that must generically be
lines (for 0 being a regular value of $u$). The parameter $\kappa$ plays a crucial role in the theory. Depending
on the value of $\kappa$, one expects different types of ``behaviors'' of vortices for minimizing configurations
relative to various constraints (boundary datas, topological constraints...) (see \cite{JT} and \cite{Ri3} for a
survey on these questions). In particular for large $\kappa$, in the so called {\it strongly repulsive limit}, the
vortices of minimizing configurations tend to minimize their length (the first order of the energy is
$2\pi\log\kappa$ times the length of the vortices) and  repel one another, which is not the case for small
$\kappa$.
 The magnetic field $dA$ plays no role in this mechanism ; to simplify the presentation
we henforth simply impose that $A=0$. We are then looking at the functional
\be
\label{7}
E_\ep(u)=\int_{B^3}|\nabla u|^2+\frac{1}{2\ep^2}(1-|u|^2)^2=\int_{B_3} e_\ep(u)
\ee
in the strongly repulsive limit $\ep\rightarrow 0$ whose critical points satisfy the non longer gauge
 invariant Ginzburg-Landau equation
\be \label{8} \Delta u+\frac{u}{\ep^2}(1-|u|^2)=0\quad\quad\mbox{
in }{\mathcal D}'({\R}^3). \ee Because we are interested in
observing finite length vortices, in view of the above remark,
 we have to restrict to critical
points satisfying $E_\ep(u_\ep)=O(\log\frac{1}{\ep})$. Equation
(\ref{8}) gives, in particular, that $\Delta u$ is parallel to
$u$. Therefore, assuming that $|u|>\frac{1}{2}$ in a ball
$B_r(x_0)$ and writing $u=|u|e^{i\phi}$, we deduce the nice
scalar elliptic equation \be \label{9}
div(|u|^2\,\nabla\phi)=0\quad\quad\mbox{ in }{\mathcal
D}'(B_r(x_0)). \ee Assuming $|u_\ep|>1/2$ on the whole $B^3$, the
compactness of solutions to (\ref{8}), in $W^{1,2}$ say, is
reduced to the compactness of the boundary data. It is then clear
that, in general, the study of the compactness of solutions
$u_\ep$ to (\ref{8}), in the limit $\ep\rightarrow 0$, involves
the study of the compactness of the sets $V_\ep=\{x\ ;\
|u_\ep(x)|< 1/2\ \}$ as well as a control of the degree of
 $u/|u|$ on arbitrary closed curves in the complement of $V_\ep$ approximating the vortices.
 (The number $1/2$ may be of course replaced by an arbitrary number between 0 and 1).

\noindent {\bf The approximate vortices {\boldmath $T_\ep$}}

In their study of the 2-dimensional version of the present problem ($B^3$ replaced by $B^2$),
 F. Bethuel, H. Brezis and F. H\'elein in \cite{BBH}, established that for solutions to
(\ref{8})
satisfying a fixed boundary condition from $\p B^2$ into $S^1$ of non-zero degree, $V_\ep$
can be covered by a uniformly bounded number of balls of radius $O(\ep)$ around which
$u/|u|$ has a uniformly bounded topological degree.
Taking then the distribution $T_\ep$ given by the sum of the  Dirac masses at the center
of the balls with the multiplicity given by the surrounding degrees of $u/|u|$, because
of the uniform bound of the mass of $T_\ep$, we may extract a converging
subsequence of atomic measures $T_{\ep'}$. The compactness of the maps $u_\ep$
is a consequence of the convergence of $T_{\ep'}$ via a classical elliptic argument.

Going back to 3 dimensions, the idea remains the same : we introduce an approximation
of the vortices and try to prove compactness. More precisely, we consider a minimal
 1-dimensional integer current $T_\ep$  in the homology class of $H_1(V_\ep,\p B^3,$ ${\Z})$ given
by the pre-image by $u_\ep$ of a regular point in $B^2_{\frac{1}{4}}$.
 Then, as in 2-dimensions, the key to the compactness of the
critical points of $E_\ep$ is now the compactness of the familly of 1-dimensional
rectifiable current $T_\ep$ constructed above.

 The energy quantization result exposed in the following subsection shows the compactness
 of the $T_\ep$.

\subsection{Quantization for G-L vortices : the {\boldmath $\eta$}-compactness}
\vskip -5mm \hspace{5mm}

The following energy quantization result introduced by the author in \cite{Ri2}
says that there exists an absolute constant $\eta>0$ such that if the
Ginzburg-Landau energy of a critical point in a given ball, suitabely renormalized, lies below this number,
 then there is no vortex passing through the ball of half radius. Precisely we have.

{\bf Theorem 5.} \cite{LR1}
\it
There exists a positive number $\eta$, such that for $\ep$ small enough
and for any critical point
of $E_\ep$ in $B^3_1$ satisfying  $|u|\le 1$, the following holds.
Let $B_r(x_0)\in B^3_{\frac{1}{2}}$, $r>\ep$ then
\[
\frac{1}{r}\int_{B_r(x_0)}e_\ep(u_\ep)\le \eta\ \log\frac{r}{\ep}\quad\quad\Longrightarrow\quad |u|>\frac{1}{2}\mbox{ in } B_{\frac{r}{2}}(x_0)\quad .
\]
\rm

The above result is optimal in the following sense : there are critical points and arbitrary $r>\ep$ such that
 $\frac{1}{r}\int_{B_r(x_0)}e_\ep(u_\ep)\le 2\pi\ \log\frac{r}{\ep}$ and $u(x_0)=0$. Such a result is reminiscent
of the $\ep-$regularity result, although one difference here is that one has to handle blowing up energy $E_\ep$ which
is not a-priori bounded as $\ep$ tends to 0.

The $\eta-$compactness was introduced first for minimizers in 3-D in \cite{Ri2} and then obtained for any
 critical points, still in 3-D, in \cite{LR1}.
Regarding dimension 2, the idea of the $\eta-$compactness is implicit in the works before \cite{Ri2}.
Making this idea explicit in \cite{Ri5} helped to substantialy simplify the existing proofs in the 2-D case.

It is a striking fact, whose explanation is beyond the scope of this
lecture, that the original proof of theorem in \cite{LR1} is based on
the same techniques used by the authors  to prove the {\it energy
identity} (theorem 3) for harmonic maps. The key of the proof is to
obtain a bound independent of $\ep$ for the density of energy
 $\frac{1}{\rho}\int_{B_\rho(x_0)}e_\ep(u)$
(for some $\rho$ between $\ep$ and $r$). The main obstacle in establishing such a bound is to control,
independently of $\ep$, the part of
the energy in $B_\rho(x_0)\setminus V_\ep$ coming from the vortices $T_\ep$ in $B_\rho(x_0)$. This is obtained by controlling,
on a generic {\it good slice} $\p B_{\rho_1}(x_0)\setminus V_\ep$ ($\rho/2<\rho_1<\rho$), the $L^{2,\infty}$ and the $L^{2,1}$
norms of $\nabla\psi$, where $\psi$ is the {\it vortex potential}  $\psi:=\Delta^{-1}T_\ep$.
Recently in \cite{BBO} F.Bethuel, H.Brezis and G.Orlandi made a nice observation by deriving, from the monotonicity
formula associated to the elliptic variational problem $E_\ep$, an $L^\infty$ bound in terms of the local density of
 energy $\frac{1}{\rho}\int_{B_\rho(x_0)}e_\ep(u)$  for the {\it vortex potential} $\psi$
in $B_\rho(x_0)\setminus V_\ep$. This bound enabled them to replace in our original proof the use of the $L^{2,\infty}-L^{2,1}$
 duality for $\nabla\psi\cdot\nabla\psi$ on a slice by a
more standard, in calculus of variation,  $L^\infty-L^1$ duality for $\psi\cdot\Delta\psi$ on $B_\rho(x_0)$ itself.
As explained in \cite{BBO}, this modification is a real improvement in the high dimensional case since it allows to end
 the
proof of the $\eta-$compactness  without having to go through an argument based on a {\it good slice} extraction
that would have been certainly more painful for $n\ge 4$. This modification of the duality used to control the energy
of the {\it vortex potential} also allowed C. Wang recently to
extend the $\eta-$compactness result from \cite{LR3} for the Ginzburg-Landau heat flow equation from 3 to 4 dimensions.

\noindent {\bf The {\boldmath $\eta$}-compactness is a compactness
result}

Considering, as mentioned above, critical points $u_\ep$ to $E_\ep$ satisfying $|u|\le 1$ and whose energy are of the order
 of $\log\frac{1}{\ep}$
($E_\ep(u_\ep)=O(\log\frac{1}{\ep}))$, a Besicovitch covering argument (see \cite{LR1}) combined with the
$\eta-$compactness above gives
\[
V_\ep\subset \cup_{j=1}^{N_\ep} B_\ep(x_j)\quad\quad\mbox{ with }\quad N_\ep=O\lf(\frac{1}{\ep}\rg)\quad .
\]
This estimate, combined with the $L^\infty$ bound of $\nabla u_\ep$ deduced from the equation and the assumption
 $|u|\le 1$,
yields the following fact:
$M(T_\ep)=O(1)$, (i.e. the mass of the approximated vortices is uniformly bounded). Combining this estimate
with the fact that, by definition, $\p T_\ep\lfloor B^3=0$, as a direct application of Federer-Fleming compactness
theorem, we deduce the compactness of the approximated vortices $T_\ep$ in the space of integer rectifiable 1-dimensional currents.

Thus, the $\eta-$compactness is the compactness of the approximated vortices $T_\ep$, which are \underbar{rectifiable currents} rather than a compactness
of \underbar{maps} $u_\ep$. The latter compactness, in $W^{1,p}$ for $p<3/2$, nevertheless  automatically follows from the previous
 one by the means of classical elliptic
arguments for a suitable class of boundary conditions for $u_\ep$ (see \cite{LR1} page 219). This class
has been recently extended to arbitrary $H^\frac{1}{2}(\p B^3,S^1)$ conditions independent of $\ep$ in \cite{BBBO}.
This idea of getting convergence of maps going through convergence of currents is in the spirit of the theory of cartesian currents
developed by M. Giaquinta, G. Modica and J. Soucek in \cite{GMS}.

\label{lastpage}

\end{document}